\documentclass[11pt]{amsart}

\newtheorem{theorem}{Theorem}[section]
\newtheorem{lem}[theorem]{Lemma}
\newtheorem{prop}[theorem]{Proposition}
\newtheorem{cor}[theorem]{Corollary}

\theoremstyle{definition}

\newtheorem{example}[theorem]{Example}

\theoremstyle{remark}

\numberwithin{equation}{section}

\newcommand{\Hor}{{\mathcal{H}}}
\newcommand{\V}{{\mathcal{V}}}
\newcommand{\ra}{\rightarrow}

\newcommand{\lb}{\langle}
\newcommand{\rb}{\rangle}

\newcommand{\mg}{\mathfrak{g}}
\newcommand{\mh}{\mathfrak{h}}

\newcommand{\bs}{\backslash}

\newcommand{\R}{\mathbb{R}}

\begin{document}

\newcommand{\spacing}[1]{\renewcommand{\baselinestretch}{#1}\large\normalsize}
\spacing{1.14}

\title{Flats in Riemannian submersions from Lie groups}

\author{Kristopher Tapp}
\address{Department of Mathematics\\ Williams College\\
Williamstown, MA 01267}
\email{ktapp@williams.edu}
\date{\today}

\begin{abstract}
We prove that any base space of Riemannian submersion from a compact Lie group (with bi-invariant metric) must have a basic property previously known for normal biquotients; namely, any zero-curvature plane exponentiates to a flat.
\end{abstract}
\maketitle
\section{Introduction}
With a few exceptions, all known examples of nonnegatively curved compact Riemannian manifolds are constructed as base spaces of Riemannian submersions from compact Lie groups with bi-invariant metrics.  The geometry of any such base space is restricted by our main theorem.  We always use $G$ to denote a compact Lie group with a bi-invariant metric.
\begin{theorem}\label{main} If $\pi:G\ra B$ is a Riemannian submersion, then
\begin{enumerate}
\item Every horizontal zero-curvature plane in $G$ projects to a zero-curvature plane in $B$.
\item Every zero-curvature plane in $B$ exponentiates to a flat (meaning a totally geodesic immersion of $\R^2$ with flat metric).
\end{enumerate}
\end{theorem}

Wilking proved Theorem~\ref{main} for biquotient submersions~\cite{Wilking}.  K. grove asked whether there exist Riemannian submersions from $G$ other than biquotient submersions~\cite{Grove}.  Our theorem forces other submersions to share a basic property of biquotient submersions, and perhaps thereby provides evidence that there are no others.  In fact, some of the ideas of our proof generalize those used by Munteanu to show that all Riemannian submersions from $G$ with one-dimensional fibers are homogeneous~\cite{Mun}.  Our generalizations could help address Grove's question with arbitrary fiber dimension.

An immediate consequence of Theorem~\ref{main} is:
\begin{cor} If $G\stackrel{\pi}{\ra} M\stackrel{f}{\ra} B$ are Riemannian submersions, then any horizontal zero-curvature plane in $M$ projects to a zero-curvature plane in $B$.
\end{cor}
In particular, Theorem~\ref{main} is true with $G$ replaced by a symmetric space, or more generally a normal biquotient.  
\section{Biquotient submersions}
In this section, we review Wilking's proof of Theorem~\ref{main} in the case of biquotient submersions, which are defined as follows.  Any subgroup $H\subset G\times G$ acts on $G$ as $(h_1,h_2)\star g=h_1gh_2^{-1}$. If this action is free, then the quotient is called a \emph{biquotient}, is denoted $G//H$, and inherits a natural metric (called a \emph{normal biquotient metric}) under which the projection $\pi:G\ra G//H$ is a Riemannian submersion (called a \emph{biquotient submersion}).  When $H\subset G\times G$ is a product, meaning $H=H_1\times H_2$, then $G//H$ is denoted $H_1\bs G/H_2$ and is called a 2-sided biquotient.

\begin{proof}[Wilking's proof for 2-sided biquotients]
Let $B=H_1\bs G/H_2$ be a 2-sided biquotient, and let $\pi:G\ra B$ be the projection map, which is a Riemannian submersion.  Let $\mg,\mh_1,\mh_2$ denote the Lie algebras of $G$, $H_1$ and $H_2$.  Let $\sigma=\text{span}\{X,Y\}\subset\mg$ be a horizontal zero-curvature plane at the identity $e\in G$.  Zero-curvature means that $[X,Y]=0$, and horizontal means that $X,Y\perp\V_e=\mh_1\oplus\mh_2$.

$F:=\exp(\sigma)\subset G$ is a flat in $G$, which we will show is everywhere horizontal.  Let $\overline{X}, \overline{Y}$ denote the left-invariant extensions of $X,Y$ to points of $F$, which are also right-invariant because $F$ is abelian. At any point $g\in F$, the vertical space is $\V_g=dR_g(\mh_1)\oplus dL_g(\mh_2)$.  The fields $\overline{X},\overline{Y}$ are orthogonal to the first factor because they are right-invariant, and to the second because they are left-invariant.  Thus, $F$ is everywhere horizontal.  Since all geodesics in $F$ are horizontal, so they project to geodesics in $\pi(F)$, which shows that $\pi(F)$ is totally geodesic, and is therefore a flat in $B$.  

The same argument works with the identity $e$ replaced by an arbitrary element of $G$, so all horizontal zero-curvatures planes in $G$ exponentiate to horizontal flats, which project to flats in $B$.  Parts (1) and (2) of the theorem follow immediately.  
\end{proof}
The above proof establishes Theorem~\ref{main} for 2-sided biquotient submersions.  Eschenburg observed that any biquotient $G//H$ is diffeomorphic to a 2-sided biquotient:
$$G//H \cong \Delta G\bs G\times G/H.$$
Metrically, an arbitrary normal biquotient of $G$ is isometric to a 2-sided normal biquotient of $G\times G$, so part (2) of Theorem~\ref{main} holds for general biquotients~\cite{Wilking}.  Part (1) for general biquotients does not seem to follow from such arguments, but does follow from Eschenburg's description of the $A$-tensor of a biquotient submersion~\cite{Esch}.
\section{Holonomy Jacobi fields}
This section contains background on Jacobi fields of $G$, particularly holonomy Jacobi fields.  We begin by reviewing the definition of holonomy Jacobi fields for a general Riemannian submersion $\pi:M\ra B$.  Let $\V,\Hor$ denote the vertical and horizontal distributions of $\pi$.  Let $A,T$ denote the fundamental tensors of $\pi$, defined as in~\cite{Besse}.

For any piecewise smooth path $\alpha(t)$ in $B$, say from $\alpha(0)=p$ to $\alpha(1)=q$,
there is a naturally associated diffeomorphism
$h_\alpha$ between the fibers $F_p:=\pi^{-1}(p)$ and $F_q:=\pi^{-1}(q)$.
This diffeomorphism maps $x\in F_p$ to the terminal point of the lift of $\alpha$ to a horizontal path in $M$ beginning at $x$.  We will call $h_\alpha$ the \emph{holonomy
diffeomorphism} associated to $\alpha$.

Let $\overline{p}\in F_p$, $v\in \V_{\overline{p}}$, and let $\overline{\alpha}(t)$ denote the horizontal lift of $\alpha(t)$ beginning at $\overline{p}$.  Denote $\overline{q}:=\overline{\alpha}(1)\in F_q$.  Notice that $d(h_\alpha)_{\overline{p}}(v)\in\V_{\overline{q}}$ equals $J(1)$, where $J(t)$ is the variational vector field along $\overline{\alpha}(t)$ obtained by lifting $\alpha$ to points of a path in $F_p$ through $\overline{p}$ in the direction of $v$.  Notice that $J(t)$ is everywhere vertical.  It is well known that for all $t\in(0,1)$:
\begin{equation}\label{JF}
J'(t) = A(\overline{\alpha}'(t),J(t))+T(J(t),\overline{\alpha}'(t)).
\end{equation}
Notice that the $A$-term in Equation~\ref{JF} is horizontal, while the $T$-term is vertical.  If $\alpha$ is a geodesic segment, then so is $\overline{\alpha}$, and the variational field $J(t)$ is a Jacobi field, which we will call a \emph{holonomy Jacobi field}.  

In summary, every vertical vector $v\in\V_{\overline{p}}$ extends to a unique holonomy Jacobi field $J(t)$ along the horizontal geodesic $\overline{\alpha}(t)$.  A basis of $\V_{\overline{p}}$ extends to a basis of holonomy Jacobi fields, which spans $\V_{\overline{\alpha}(t)}$ for all $t$, and thus determine the vertical distribution along $\overline{\alpha}$.  On the other hand, since Jacobi fields are determined by their initial values and initial derivatives, Equation~\ref{JF} shows that the $A$ and $T$ tensors at the single point $\overline{p}$ determine the vertical distribution along all horizontal geodesics emanating from $\overline{p}$.  This provides substantial rigidity, especially in situations where the Jacobi equation can be solved.

Next, we review the solution to the Jacobi equation on a compact Lie group $G$ with a bi-invariant metric.  Let $\mg$ denote the Lie algebra of $G$, and let $X\in\mg$ be a unit-length vector.  Let $\gamma(t)=\exp(tX)$ be the geodesic in the direction of $X$.  Following~\cite{Mun}, we will characterize the Jacobi fields along $\gamma(t)$.

Decompose $\mg$ into eigenspaces of $(\text{ad}_X)^2$:
\begin{equation}\label{root}\mg = V_0 + \sum_{i=1}^{l} V_i,\end{equation}
with corresponding non-positive eigenvalues $0=\lambda_0>\lambda_1>\cdots>\lambda_l$.  Notice that for any $v_i\in V_i$, the sectional curvature of $\text{span}\{X,v_i\}$ is $k_i:=-\lambda_i/4$.  Notice that if $X$ is regular, then $V_0$ is the maximal abelian subalgebra containing $X$, and Equation~\ref{root} is the root space decomposition.

Jacobi fields along $\gamma(t)$ all have the form:
\begin{equation}\label{J}
J(t) = E_0 + t\cdot F_0 + \sum_{i=1}^{l}\left(\cos(\sqrt{k_i}t)\cdot E_i + \sin(\sqrt{k_i}t)\cdot F_i\right)
\end{equation}
where $E_i,F_i\in V_i$ are parallel translated along $\gamma(t)$.  In fact, the above equation describes the unique Jacobi field along $\gamma(t)$ with initial data:
$$J(0) = E_0 + \sum_{i=1}^{l}E_i,\,\,\,\,\,\,\,\,
  J'(0) = F_0 + \sum_{i=1}^{l}\sqrt{k_i}\cdot F_i.$$

\section{Holonomy Jacobi fields are bounded}
The crux of the proof of Theorem~\ref{main} is:
\begin{lem}
For a Riemannian submersion $\pi:G\ra B$, every holonomy Jacobi field remains bounded in norm.
\end{lem}
This was proven in~\cite{Mun} assuming that the fibers are one-dimensional.  For general Riemannian submersions between compact manifolds, holonomy Jacobi fields may grow unbounded.  Examples of this phenomenon will be explored in Section 6.
\begin{proof}
Let $X\in\mg$ be a horizontal vector, and let $\gamma(t)=\exp(tX)$ be the geodesic in the direction of $X$.  We wish to prove that holonomy Jacobi fields along $\gamma$ remain bounded.

Define $\Omega$ as the set of vertical vectors at $e$ which determine parallel holonomy Jacobi fields along $\gamma$:
$$\Omega:=\{v\in\V_e\mid \text{the holonomy Jacobi field }J(t)\text{ along }\gamma(t)\text{ with } J(0)=v \text{ is parallel}\}.$$
Notice that $\Omega$ is a (possibly trivial) subspace of $\V_e$, and every vector in $\Omega$ commute with $X$.

If $\gamma$ is a closed geodesic, it is easy to show that $\Omega$ parallel transports to itself along $\gamma$.  We next prove that, in a certain limit sense, $\Omega$ returns to itself even when $\gamma$ is not closed.  The closure of the image of $\gamma$ is a torus in $G$.  Choose any sequence $\{t_i\}$, with $t_i\ra\infty$, such that $\gamma(t_i)\ra e$, and $\gamma'(t_i)\ra X$.  Let $V(t)$ be the parallel transport along $\gamma(t)$ of a vector $V(0)\in\Omega$.  Consider the sequence $\{V(t_i)\}$ in $TG$.  Let $\tilde{V}$ denote any sub-limit of this sequence. By continuity of the vertical distribution, $\tilde{V}\in\V_e$.  We claim that $\tilde{V}\in\Omega$.  First notice that $[X,\tilde{V}]=0$ by continuity, since for each $i$, $[dL_{\gamma(t_i)}\gamma'(t_i),dL_{\gamma(t_i)}V(t_i)]=0$.  So by Equations~\ref{JF} and ~\ref{J}, it remains to verify that $A(X,\tilde{V})+T(\tilde{V},X)=0$, which is a consequence of the continuity of $A$ and $T$, since $A(\gamma'(t_i),V(t_i))+T(V(t_i),\gamma'(t_i))=0$ for each $i$.

The above closure property implies a sense in which $\Omega$ parallel translates along $\gamma(t)$ to itself at $t=\infty$.  Indeed, take an orthonormal basis $\{V_k\}$ of $\Omega$. After restricting to a subsequence of $\{t_i\}$, we can insure that, for each basis element $V_k$, parallel translated along $\gamma(t)$, the sequence $V_k(t_i)$ converges to some limit $\tilde{V}_k\in\V_e$.  These limits $\{\tilde{V}_k\}$ are orthonormal, and by the above closure property they lie in $\Omega$, so they form an orthonormal basis of $\Omega$.  Henceforth, we let $\{t_i\}$ denote this subsequence, and let $\Omega(t)$ denote the parallel transport of $\Omega$ along $\gamma(t)$, so that $\Omega(t_i)\ra\Omega$.

Let $\Omega^{\perp}$ be the orthogonal compliment of $\Omega$ in $\V_e$.  If $J(t)$ is a holonomy Jacobi field along $\gamma(t)$ with $J(0)\in\Omega^\perp$, then for all $V\in\Omega$,
$$\lb J'(0),V\rb = \lb T(J(0),X),V\rb = \lb T(V,X),J(0)\rb = \lb 0,J(0)\rb = 0.$$
Thus, $J'(0)\perp\Omega$.  It follows from Equation~\ref{J} that $J(t)\perp\Omega(t)$ for all $t$.

So to finish the proof, it will suffice to show that any holonomy Jacobi fields $J(t)$ along $\gamma(t)$ with $J(0)\in\Omega^{\perp}$ remain bounded.  Let $J(t)$ be such a holonomy Jacobi field.  Assume that $|J(t)|$ is unbounded, which means that $F_0\neq 0$ in Equation~\ref{J}.  Notice that $F_0\perp\Omega$ because $J'(0)\perp\Omega$.
Let $F(t)$ denote the parallel transport of $F(0):=F_0/|F_0|$ along $\gamma(t)$.  Notice that $F(t)\perp\Omega(t)$ for all $t$.

By Equation~\ref{JF},
$$|J'(t)| = |J(t)|\cdot\left|A\left(\gamma'(t),\frac{J(t)}{|J(t)|}\right)+T\left(\frac{J(t)}{|J(t)|},\gamma'(t)\right)\right|.$$
Notice that $|J'(t)|$ is bounded from above by Equation~\ref{J}, while $|J(t)|$ grows unbounded, and in fact is asymptotic to a linear function with slope $|F_0|$.  It follows that:
$$\lim_{t\ra 0} \left|A\left(\gamma'(t),\frac{J(t)}{|J(t)|}\right)+
T\left(\frac{J(t)}{|J(t)|},\gamma'(t)\right)\right|=0.$$
For large $t$, $J(t)/|J(t)|$ is close to $F(t)$, and in particular $F(t)$ is close to vertical, so we can re-write this as:
\begin{equation}\label{lim}\lim_{t\ra 0} |A(\gamma'(t),F(t)^\V)+T(F(t)^\V,\gamma'(t))|=0.
\end{equation}

After passing to a subsequence, we can insure that $F(t_i)$ converges to some limit $\tilde{F}\in\V_e$.  By continuity, $[X,\tilde{F}]=0$.  Finally, Equation~\ref{lim} together with continuity imply that:
$$|A(X,\tilde{F})+T(\tilde{F},X)|=0.$$
Equation~\ref{JF} implies that the holonomy Jacobi field along $\gamma(t)$ beginning at $\tilde{F}$ has zero initial derivative, and is therefore parallel, so $\tilde{F}\in\Omega$.

On the other hand, $\tilde{F}\in\Omega^\perp$, since for each $i$, $F(t_i)\perp\Omega(t_i)$, and we saw previously that $\Omega(t_i)\ra\Omega$.  This provides a contradiction.
\end{proof}
\section{proof of main theorem}
\begin{proof}[Proof of Theorem~\ref{main}]
Let $\pi:G\ra B$ be a Riemannian submersion.  Let $\sigma=\text{span}\{X,Y\}$ be a horizontal zero-curvature plane at $e$.  Zero-curvature means that $[X,Y]=0$.  The image $F:=\exp(\sigma)$ is a flat in $G$, which we wish to prove is everywhere horizontal.  Let $\gamma(t)=\exp(tX)$, and let $Y(t)$ be the parallel transport of $Y$ along $\gamma(t)$.  A basis for $T_{\gamma(t)} F$ is $\{\gamma'(t),Y(t)\}$.  Since $\gamma'(t)$ is clearly horizontal, it suffices to prove that $Y(t)$ is horizontal.  Equivalently, we must show that $Y(t)$ is orthogonal to any holonomy Jacobi field $J(t)$ along $\gamma(t)$.  By Equation~\ref{J}, the holonomy Jacobi field $J(t)$ remains orthogonal to $Y(t)$ if and only if $J'(0)\perp Y$.  But if $J'(0)$ were not orthogonal to $Y$, then $J(t)$ would grow unbounded, contradicting the previous lemma.

The same argument works with the identity $e$ replaced by an arbitrary element of $G$, so all horizontal zero-curvatures planes in $G$ exponentiate to horizontal flats, which project to flats in $B$.  The theorem follows.
\end{proof}
\section{Examples of unbounded holonomy Jacobi fields}
In this section, we mention conditions under which holonomy Jacobi fields remain bounded for an arbitrary Riemannian submersion $\pi:M\ra B$.  We first show an example where they do not, even though $M$ and $B$ are compact.

\begin{example}\label{E}
Let $B$ and $F$ be compact Riemannian manifold.  Endow $B\times F$ with the product metric, and let $\pi:B\times F\ra B$ denote the projection, which is a Riemannian submersion.  Let $X$ and $Y$ denote vector fields on $B$ and $F$ respectively.  Let $\mu(t)$ be a smooth function which equals 0 for $t\leq 0$ and is elsewhere positive.  Define the following distribution on $B\times F$:
$$\Hor(p,q) = \{(A,0)\mid A\perp X(p)\}\oplus\text{span}\left\{\left(X(p),\mu(|X(p)|)\cdot Y(q)\right)\right\}.$$
Notice that if $X(p)=0$, then $\Hor(p,q)$ agrees with the original horizontal distribution of $\pi$.  The factor $\mu(|X(p)|)$ insures that $\Hor$ is smooth at the boundary of the support of $(X,0)$.  There exists a unique metric on $B\times F$ such that $\pi$ remains a Riemannian submersion, the metric on the vertical distribution is unchanged, and $\Hor$ becomes the new horizontal distribution.  All fibers of $\pi$ remain naturally isometrically identified with $F$.  The holonomy diffeomorphism associated to a path $\gamma:[a,b]\ra B$ is the map from $F$ to $F$ which flows along $Y$ for time $T=\int_a^b\overline{\mu}(\lb \gamma'(t),X(\gamma(t))\rb)$, where $\overline{\mu}$ is the odd function that agrees with $\mu$ for positive $t$.  It is easy to choose $F$ and $Y$ such that the Lipschitz constant for this flow will grow unbounded as $T\ra\infty$.  If $\gamma:[0,\infty)\ra B$ is a geodesic such that $\int_0^\infty\overline{\mu}(\lb \gamma'(t),X(\gamma(t))\rb)=\infty$, then there must exist a holonomy Jacobi fields along the lift of $\gamma(t)$ which grown unbounded.
\end{example}

If the fibers of a Riemannian submersion are totally geodesic, then any holonomy Jacobi field has constant length, and therefore has bounded norm.  A more general sufficient condition involves the notion of compact holonomy.  Fix $p\in B$ and define the \emph{holonomy group of $\pi$} as the group, $\Phi$,
of diffeomorphisms of the fiber $F_p=\pi^{-1}(p)$ which occur as $h_\alpha$ for a piecewise smooth loop $\alpha$ in $B$ at $p$.  $\Phi$ is clearly independent of the choice of $p$, up to group isomorphism, when $M$ is connected.  $\Phi$ is not necessarily a finite dimensional Lie group.  We say that $\pi$ has \emph{compact holonomy} if $\Phi$ is a compact finite dimensional Lie group.  The following is proven in~\cite{T}:
\begin{prop}
If $\Phi$ is compact, then there exists a constant $L$ such that all holonomy diffeomorphisms of $\pi$ satisfy the Lipschitz constant $L$.  In particular, any holonomy Jacobi field remains bounded.
\end{prop}

In Example~\ref{E}, the holonomy group is isomorphic to $\R$, which is finite dimensional but not compact.

\bibliographystyle{amsplain}

\end{document}